\documentclass[conference]{IEEEtran}

\ifCLASSINFOpdf
   \usepackage[pdftex]{graphicx}
   \DeclareGraphicsExtensions{.pdf,.jpeg,.png}
\else
\fi
\usepackage{xcolor}
\usepackage{amsthm}

\newtheorem{definition}{Definition}[section]
\usepackage{amsmath}
\usepackage{amsfonts}
\usepackage{soul}

\usepackage{url}

\usepackage[]{fancyhdr} %
\newcommand{\changefont}{\fontsize{9}{9}\selectfont}
\ifCLASSOPTIONcompsoc
  \usepackage[caption=false,font=normalsize,labelfont=sf,textfont=sf]{subfig}
\else
  \usepackage[caption=false,font=footnotesize]{subfig}
\fi
\IEEEoverridecommandlockouts
\begin{document}

%
\title{Market mechanism to enable grid-aware dispatch of Aggregators in radial distribution networks}

\author{\IEEEauthorblockN{Nawaf Nazir}
\IEEEauthorblockA{\textit{Energy and Environment Directorate}\\Pacific Northwest National Laboratory\\
Richland, USA\\
nawaf.nazir@pnnl.gov}
\and
\IEEEauthorblockN{Mads Almassalkhi}
\IEEEauthorblockA{\textit{Department of Electrical and Biomedical Engineering}\\University of Vermont\\
Burlington, USA\\
malmassa@uvm.edu} \thanks{This work was supported by the U.S. Department of Energy’s Advanced Research Projects Agency—Energy (ARPA-E) Award DE-AR0000694 and NSF Award ECCS-2047306 is gratefully acknowledged. 
}
}


%





\maketitle
\thispagestyle{fancy}
\pagestyle{fancy}


\begin{abstract}
This paper presents a market-based optimization framework wherein Aggregators can compete for nodal capacity across a distribution feeder and guarantee that allocated flexible capacity cannot cause overloads or congestion. This mechanism, thus, allows Aggregators with allocated capacity to pursue a number of services at the whole-sale market level to maximize revenue of flexible resources.
 Based on Aggregator bids of capacity (MW) and network access price (\$/MW), the distribution system operator (DSO) formulates an optimization problem that prioritizes capacity to the different Aggregators across the network while implicitly considering AC network constraints. This grid-aware allocation is obtained by incorporating a convex inner approximation into the optimization framework that prioritizes hosting capacity to different Aggregators. We adapt concepts from transmission-level capacity market clearing, utility demand charges, and Internet-like bandwidth allocation rules to distribution system operations by incorporating nodal voltage and transformer constraints into the optimization framework. Simulation based results on IEEE distribution networks showcase the effectiveness of the approach.
\end{abstract}



%
\IEEEpeerreviewmaketitle

\section{Introduction}
 The addition of large-scale renewable energy in power systems has led to increased volatility in supply that can negatively affect the reliability of the distribution grid. To counter this volatility in renewable generation, flexibility offered by demand-side resources (such as electric vehicles and controllable loads) has shown promise~\cite{mathieu2012state,desrochers2019real}. Unlike traditional MW-scale thermal generators, flexibility in medium- and low-voltage networks will be composed of thousands of kW-scale customer-owned devices that are distributed throughout the medium- and low-voltage networks and will require careful coordination and communication. Many of these customer-owned devices can be controlled as an aggregate resource by an Aggregator. Uncoordinated operation of these Aggregators can lead to reliability issues and congestion within the distribution network~\cite{ross2021strategies}. This motivates the need to coordinate Aggregators across the network in a manner that guarantees satisfaction of AC network constraints (i.e., a grid-aware allocation of Aggregator flexibility). Based on this, we cast the grid-aware allocation of Aggregators as a market-based, convex optimization problem, which metes out the effective hosting capacity of a radial feeder based on Aggregators' capacity bids (MW, \$/MW).

\begin{figure}[!t]
    \centering
    \subfloat{\includegraphics[width=0.48\linewidth]{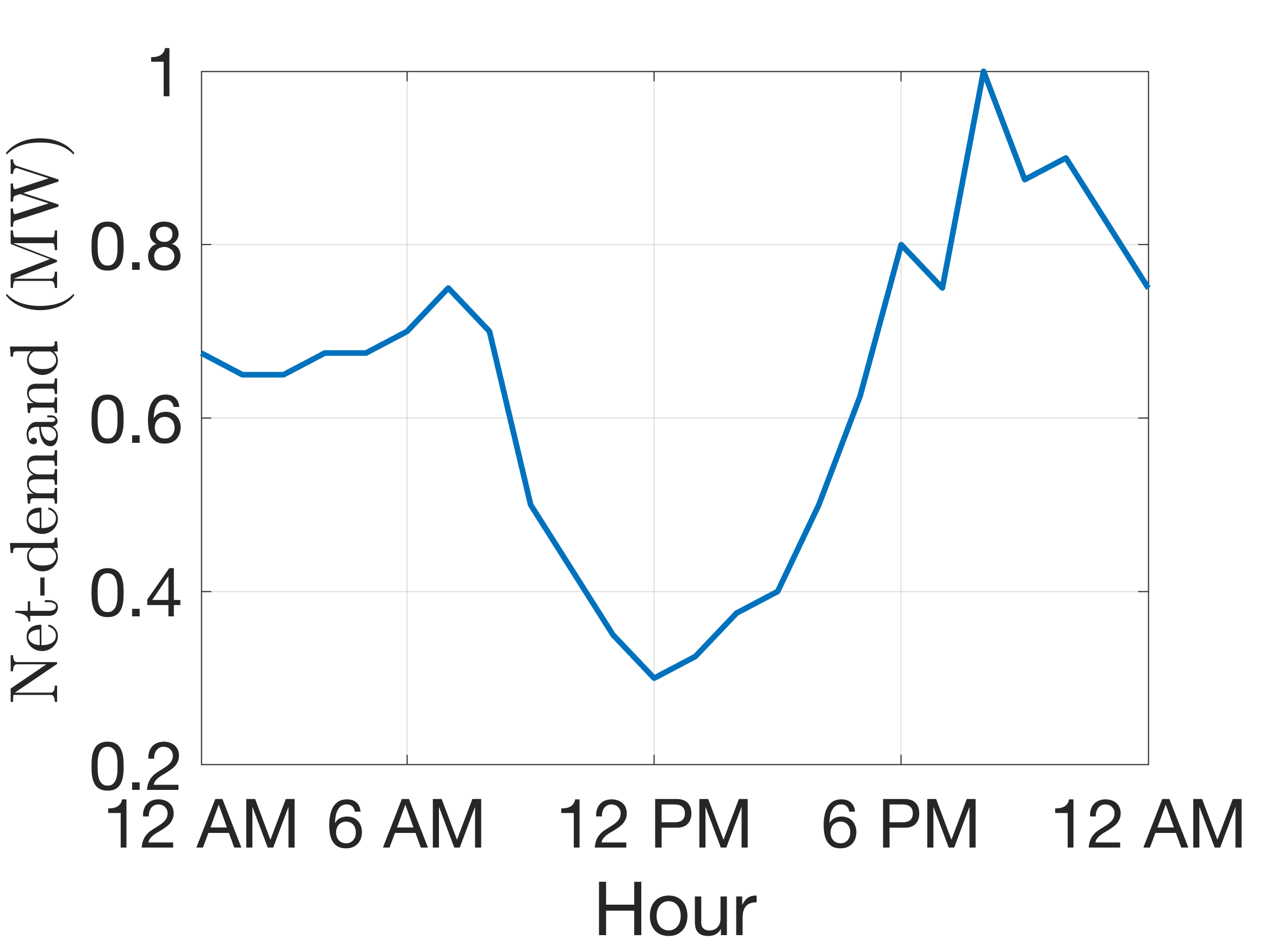}\label{fig:net_demand}}
    \subfloat{
    \includegraphics[width=0.48\linewidth]{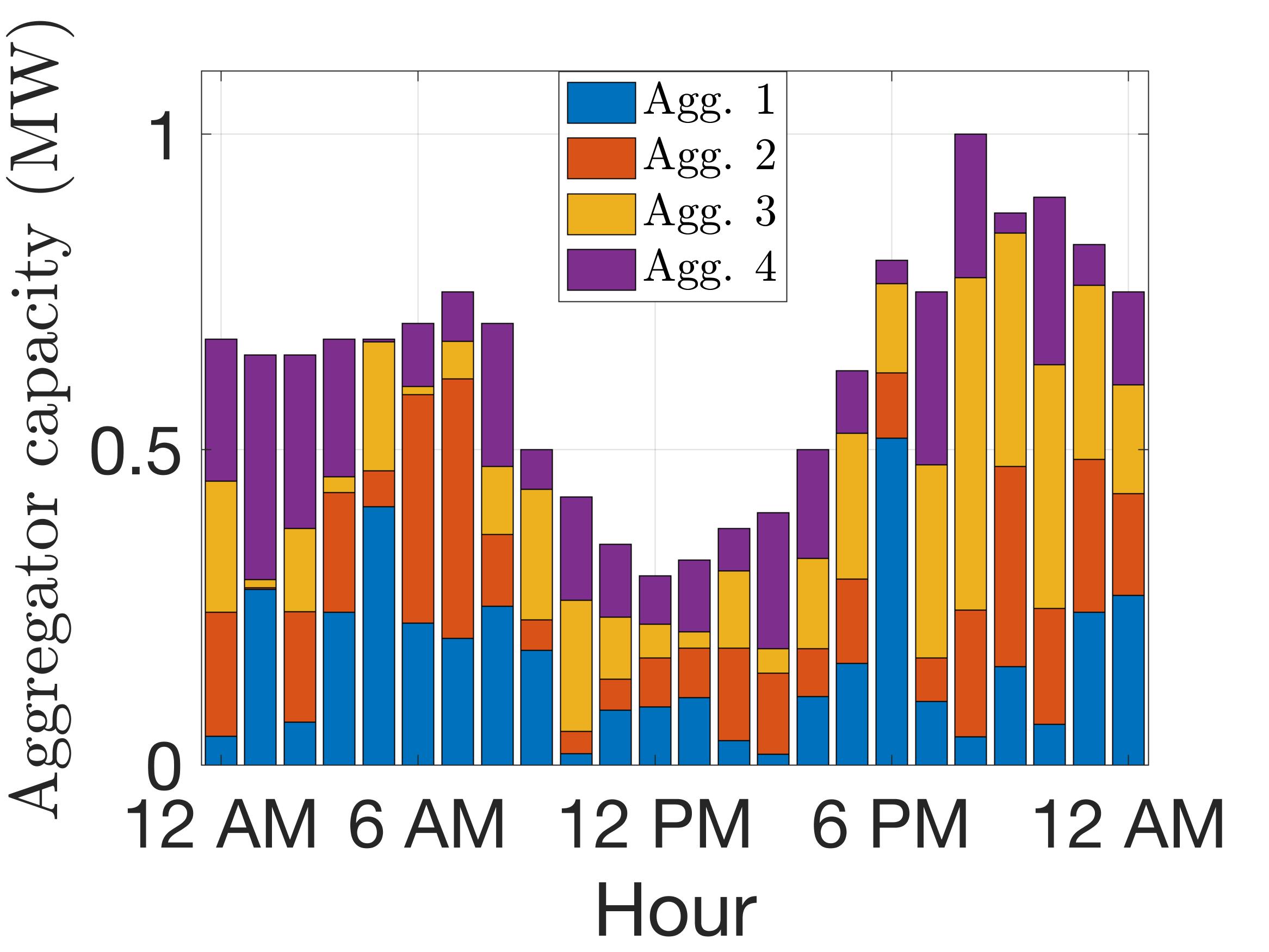}\label{fig:Agg_capacity}}
\caption{(a) Day-ahead net-demand at a single primary distribution node, (b) Day-ahead aggregator flexibility capcity to meet the net-demand.}
\end{figure}

The concept of hosting capacity (HC) is not new and has been used by utilities to site and size solar PV and EV charging station installations in a grid-aware manner. However, HC methods tend to be piecemeal and static and consider worst-case annual demand profiles, which can lead to overly conservative HC values and an under-utilization of distribution network capacity. Demand-side flexibility (e.g., from controlled behind-the-meter loads, EVs chargers, and batteries) can also be thought of as capacity, however, such flexible assets represent dynamic resources and require a re-thinking of the conventional (static) HC concept~\cite{nando_envelopes}. Thus, our aim in this work is to determine a distribution feeder's dynamic hosting capacity and how to prioritize nodal access in a grid-aware manner when multiple Aggregators offer flexibility at different locations in the network. The authors'  recent work on determining nodal hosting capacity developed a novel convex inner approximation (CIA) of the set of admissible active power injections~\cite{nazir2020gridaware} that embeds network voltage and line flow constraints within the nodal (net) injection bounds. Conventional OPF methods can also dispatch flexibility in a grid-aware manner but require full network and device data, which is only reasonable when the Distribution network operator (DNO) is also the aggregator (i.e., vertically integrated)~\cite{nazir_energies}.
In this paper, we aim to utilize the CIA formulation to determine the dynamic hosting capacity of Aggregators by formulating it as a market-clearing mechanism.

Such an approach can provide a grid-aware interface between Aggregators (that manage behind-the-meter devices) and a distribution network operator (DNO) or Utility (that is responsible for maintaining system reliability and security). Another key benefit of the presented approach is that it establishes a clear separation between Aggregators' customer-owned devices and their associated private device data and the Utility's sensitive network data.  That is, the grid-aware allocation is achieved without Aggregators needing any grid data nor Utility needing any device or personal usage data. Finally, since the DSO's role focuses on grid-aware allocation of Aggregator flexibility, it allows the Aggregators to quantify their dispatchable flexibility \textit{a priori} to maximize revenues from participating in multiple market opportunities on different timescales (e.g., ancillary services, capacity, frequency regulation, arbitrage).

\subsection{Literature survey}
Grid-aware coordination of distributed resources has previously been studied in several works in literature, such as~\cite{molzahn2019grid} that provides a certificate test whether grid-aware control is required,~\cite{ross2021strategies} that provides various strategies for network safe coordination of Aggregators,~\cite{baker_cognizant} that describes local control laws that respond to real-time changes in voltage levels and~\cite{nando_envelopes} that develops network constraint dependent import/export envelopes for Aggregators to operate under. Related to the optimal power flow based methods, several works in literature such as~\cite{nazir_3ph,nazir_energies} have developed centralized optimization based frameworks to dispatch flexible resources in three-phase networks. In this manuscript, we develop a method for a DSO to optimally allocate network capacity between competing Aggregators, such that no Aggregator dispatch will negatively impact grid reliability. Thus, the DSO's allocation of Aggregators' flexibility is grid-aware and enables Aggregators to participate with their allocated flexibility in any whole-sale TSO markets without the potential to cause network congestion in the radial distribution feeder. We achieve this by formulating the optimization problem as a day-ahead market-clearing mechanism.

Day ahead market-clearing mechanisms have been studied in the literature with~\cite{Bai_DLMP} proposing a distribution locational marginal pricing (DLMP) for active and reactive powers. These distribution-level market mechanism generally employ DLMPs within transactive energy schemes that broadcast variable DMLPs to devices to incentivize participation~\cite{Hammerstorm_transactive}. Combining grid-aware DER control with market-clearing methods has also been studied with~\cite{hu2018aggregator} proposing a network constrained transactive energy framework to coordinate the flexibility. These works, however, either neglect the distribution network constraints, or utilize a simplified linear model, that does not provide guarantees on system reliability. In this paper, we present a grid-aware market-clearing mechanism that allocates and prioritizes flexibility to competing Aggregators while avoiding network congestion.
Based on the above, this manuscript contributes to the literature as follows:
\begin{enumerate}
    \item We adapt the CIA formulation within a market mechanism to allocate and prioritize flexibility amongst competing Aggregators located across a distribution network.
    \item We introduce a two-step market mechanism where step~1 checks the feasibility of hosting the available Aggregator flexibility and in case of expected congestion, step~2 allocates flexibility to Aggregators based on their bids.
    \item We extend this mechanism to account for uncertainty in background demand with a robust formulation.
\end{enumerate}

\section{Problem definition and System Model}
\begin{figure}[!t]
\centering
\includegraphics[width=0.45\textwidth]{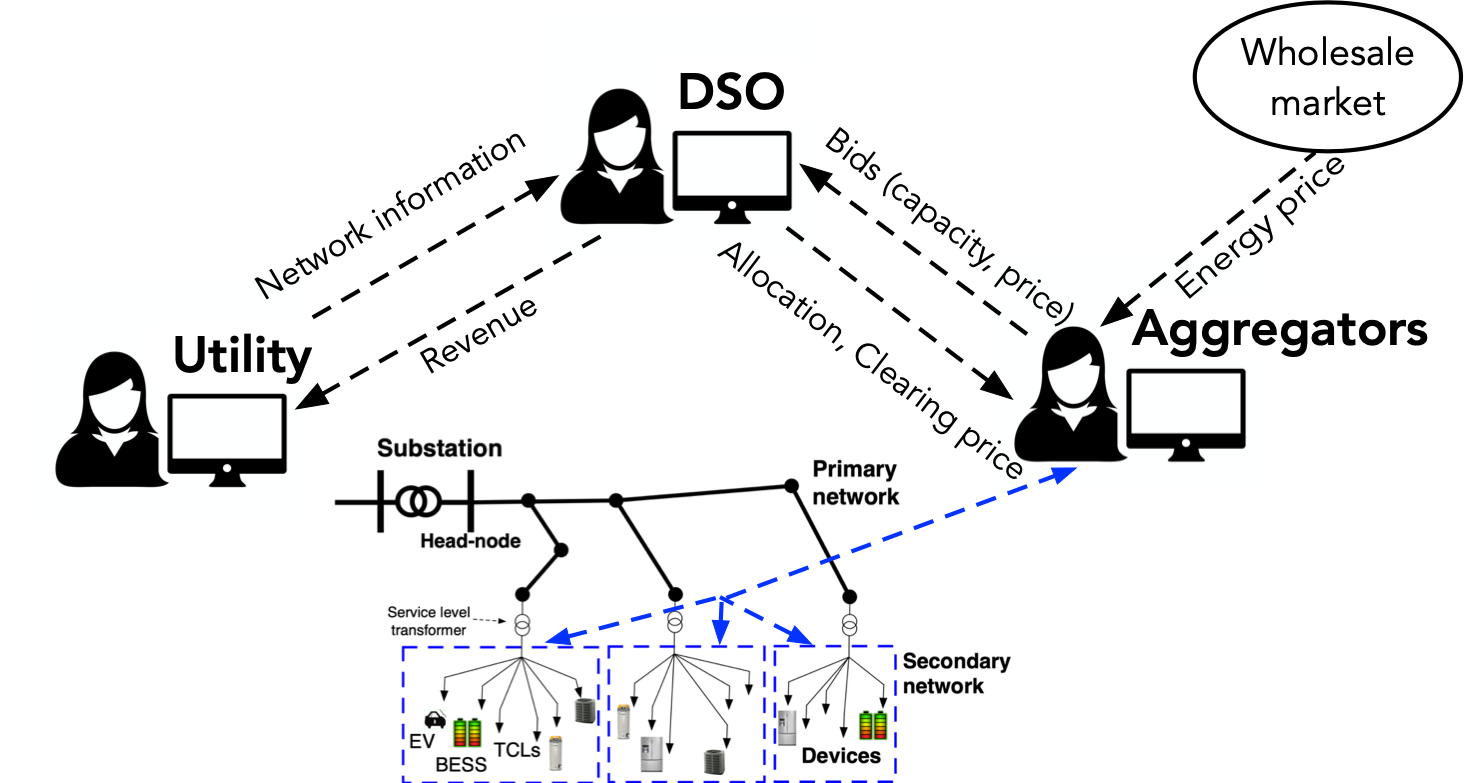}
\caption{\label{fig:market_arch} Illustration of market architecture showing information exchange between different entities.}
\end{figure}

The proposed market architecture is illustrated in Fig.~\ref{fig:market_arch}, showing the different actors within the distribution system (DSO, DNO, and Aggregators) and highlighting their information exchanges. The DSO is assumed to be an independent authority that manages Aggregators' bids (MW, \$/MW) and DNO's network data and clears the flexibility market when congestion is expected. The DNO owns and maintains the infrastructure and is compensated by the Aggregators to use the network during periods of congestion. The Aggregators pay a fee (set by the DSO's market clearing mechanism) to the DNO in order to access the distribution network capacity at their cleared allocation, \textcolor{black}{called the market-clearing price.}
\textcolor{black}{
\begin{definition}[Market-clearing price]
A uniform price auction equilibrium solution at the intersection of the distribution network capacity and the Aggregators' bids.
\end{definition}
}
The Aggregators then use their allocation to actively participate in different valuable wholesale markets. As part of the bid structure, the Aggregators provide the maximum capacity charges (\$/MW) that they are  willing to pay for access to the distribution network, which could either be a single price per feeder or nodal-level prices.
These prices are realized in the same manner as demand charges but within the context of a network's hosting capacity.
\textcolor{black}{In addition, the network capacity allocation concept is analogous to the ones utilized in wireless network bandwidth allocation, where ISPs provide Mbps capacity allocations to its users ~\cite{frolik2004qos}}.

Thus, based on internal costs of deploying demand-side flexibility, including the costs to acquire, deploy, and operate connected devices and expected revenue from delivering whole-sale market services, the Aggregators need to compute their own prices.
The problem considered in this paper is to allocate flexible hosting capacity across a distribution network to multiple Aggregators while accounting for the nonlinear AC physics and operational constraints. Specifically, we break down the problem into two steps: i) feasibility check ii) prioritization of flexibility. In Step~1, the DSO receives each Aggregators' capacity bid (MW) and the DNO's network model to check if any Aggregator actions can cause network congestion. If Step~1 is feasible, then all flexibility is accepted and Aggregators can dispatch freely while the DNO is guaranteed that Aggregator actions will not impact reliability. If Step~1 is infeasible, the DNO's network cannot host all Aggregators' flexible capacity and the DSO uses the prices (\$/MW) from the Aggregator's bid to prioritize and allocate Aggregators at each node in the network.
In the next section, we provide the mathematical model of the distribution grid that will be utilized in solving this market optimization problem.

\subsection{Nonlinear AC model}\label{sec:math_model}
\begin{figure}[!t]
\centering
\includegraphics[width=0.38\textwidth]{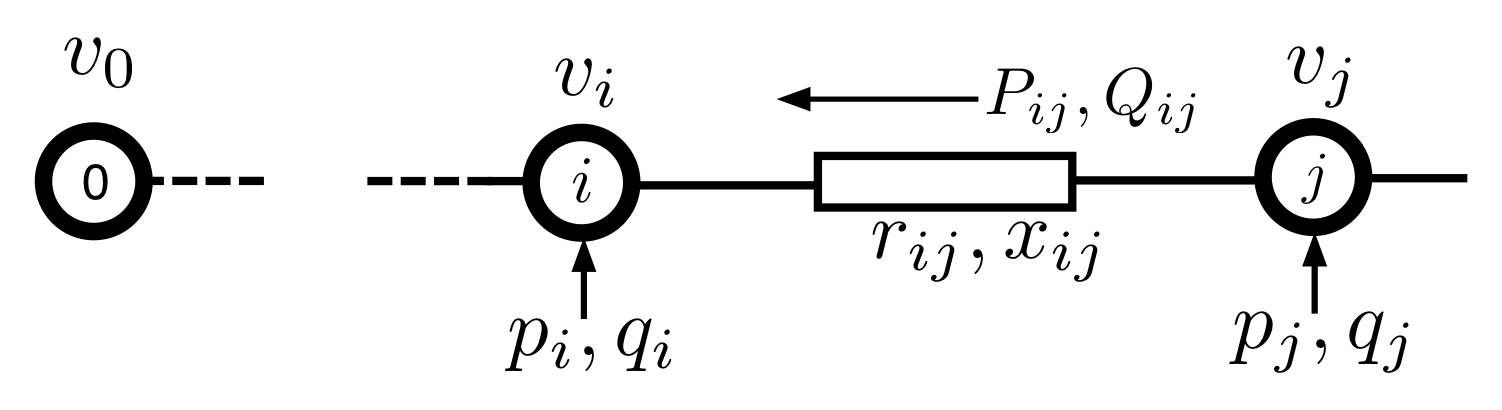}
\caption{\label{fig:radial_network} Nomenclature for a radial distribution network~\cite{heidari2017non}.}
\end{figure}

Consider a balanced, radial distribution network, shown in Fig.~\ref{fig:radial_network}, as an undirected graph $\mathcal{G}=\{\mathcal{N}\cup\{0\}, \mathcal{L}\}$ consisting of a set of $N+1$ nodes with $\mathcal{N}:=\{1, \hdots, N\}$ and a set of $N$ branches $\mathcal{L}:=\{1, \hdots, N\} \subseteq \mathcal{N}\times \mathcal{N}$, such that $(i,j)\in \mathcal{L}$, if nodes $i,j$ are connected. 
Node $0$ is assumed to be the substation node with a fixed voltage $V_0$.
If $V_i$ and $V_j$ are the voltage phasors at nodes $i$ and $j$ and $I_{ij}$ is the current phasor in branch $(i,j)\in \mathcal{L}$, then define $v_i:=|V_i|^2$, $v_j:=|V_j|^2$ and $l_{ij}:=|I_{ij}|^2$. Let $P_{ij}$ ($Q_{ij}$) be the active (reactive) power flow from node $j$ to $i$, let $p_j$ ($q_j$) be the active (reactive) power injections into node $j$, and let $r_{ij}$ ($x_{ij}$) be the resistance (reactance) of branch $(i,j)\in \mathcal{L}$, which means that the branch impedance is given by $z_{ij}:=r_{ij}+\mathbf{i}x_{ij}$. Then, for a radial network, the relation between node voltages and power flows is given by the \textit{DistFlow} equations $\forall (i,j)\in \mathcal{L}$:

\begin{subequations}\label{eq:dist_flow}
\begin{align}
v_j=&v_i+2r_{ij}P_{ij}+2x_{ij}Q_{ij}-|z_{ij}|^2l_{ij} \label{eq:volt_rel}\\
P_{ij}=&p_j+\sum_{h:h\rightarrow j}(P_{jh}-r_{jh}l_{jh}) \label{eq:real_power_rel}\\
Q_{ij}=&q_j+\sum_{h:h\rightarrow j}(Q_{jh}-x_{jh}l_{jh}) \label{eq:reac_power_rel}\\
l_{ij}(P_{ij},Q_{ij},v_j)  =& \frac{P_{ij}^2+Q_{ij}^2}{v_j}, \label{eq:curr_rel}
\end{align}
\end{subequations}
where nodal power injections are $p_j:=p_{\text{g},j}-P_{\text{L},j}$ and $q_j:=q_{\text{g},j}-Q_{\text{L},j}$ with $p_{\text{g},j}$ ($q_{\text{g},j}$) as the controllable active (reactive) flexibility and $P_{\text{L},j}$ ($Q_{\text{L},j}$) is the uncontrollable active (reactive) demand.  

We seek to determine if the network's available flexible capacity, e.g., $p_g^+ \in \mathbb{R}^N$, can be hosted by the distribution network. That is, for all $p_g \in [0,p_g^+]$, do the corresponding voltages and currents remain within limits, i.e., $v_j(p_g)\in [\underline{v_j},\overline{v_j}]\quad \forall j\in \mathcal{N}$ and $l_{ij}(p_g) \in [\underline{l_{ij}},\overline{l_{ij}}] \quad \forall (i,j)\in \mathcal{L}$. We can formulate the problem as:
{\color{black}
\begin{align}
\text{(NLP)}  \quad s_i^*=\arg\min_{s_i^+} \  \sum_{i=1}^{N}s_i^+&\\
\text{s.t.}  \quad  \eqref{eq:volt_rel}-\eqref{eq:curr_rel}\\
    p_i=p_{\text{g},i}^+-P_{\text{L},i}-s_i^+ & \quad \forall i\in \mathcal{N} \\
    \underline{v_i}\le v_i \le \overline{v_i} & \quad \forall i\in \mathcal{N}   \label{eq:P1_V}\\
    \underline{l_{ij}}\le l_{ij} \le \overline{l_{ij}} & \quad \forall (i,j)\in \mathcal{L} \label{eq:P1_lmax}\\
    s_i^+ \ge 0 & \quad \forall i\in \mathcal{N}
\end{align}
}
\textcolor{black}{where $p_{\text{g},i}^+$ is the available upper flexibility at node $i$ and $s_i^+$ is the slack variable when considering upper flexibility at node $i$. A similar optimization problem can be formulated to determine the slack variable $s_i^-$ when considering the lower flexibility $p_{\text{g},i}^-$.} In either case, if the flexibility $p_{\text{g},i}$ can be accommodated by the network without violating any constraints, then $s_i^*=0 \quad \forall i\in \mathcal{N}$, and there is no possibility of network congestion and a market mechanism is not required to prioritize flexibility in the network and hence the DSO does not have to worry about network reliability. If however, $\exists i \in \mathcal{N}$, s.t., $s_i^*>0$, then uncoordinated dispatch of flexible resources can lead to network congestion and violation of network constraints. In order to avoid such a situation, we need to prioritize the available flexibility.

However, solving (NLP) in general represents a technical challenge, as the \textit{DistFlow} equations governing the distribution system are non-linear due to~\eqref{eq:curr_rel}, making the optimization problem (NLP) non-convex. In order to overcome this challenge, many works in literature have considered approximations and convex relaxation techniques~\cite{molzahn2019survey}. Amongst the approximation methods is the so-called linear \textit{DistFlow} model or \textit{LinDist}~\cite{baran1989optimal}, whereas the convex relaxation techniques include a host of different methods~\cite{taylor2015convex}. In the next section, we will discuss the shortcomings of these methods, particularly as it relates to the original (NLP), and then provide an alternative method that overcomes those shortcomings.


\subsection{Linear and relaxation based models}
This section presents the effects of using linear and convex relaxation based models when determining maximum active power injections. First, consider a linear model.

\subsubsection{Shortcomings of \textit{LinDist}}
\begin{figure}[!t]
    \centering
    \subfloat{
    \includegraphics[width=0.48\linewidth]{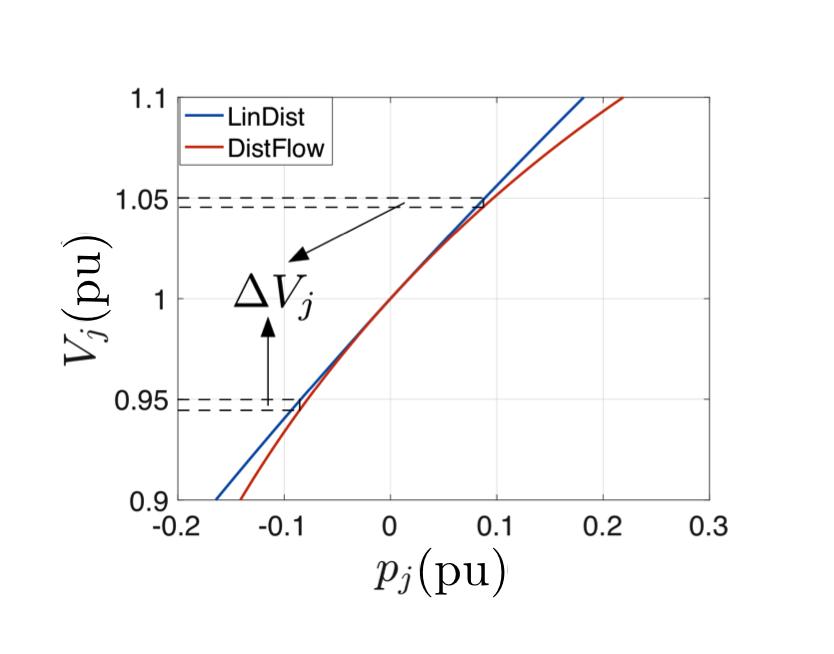}\label{fig:QvsV}}
    \subfloat{
    \includegraphics[width=0.48\linewidth]{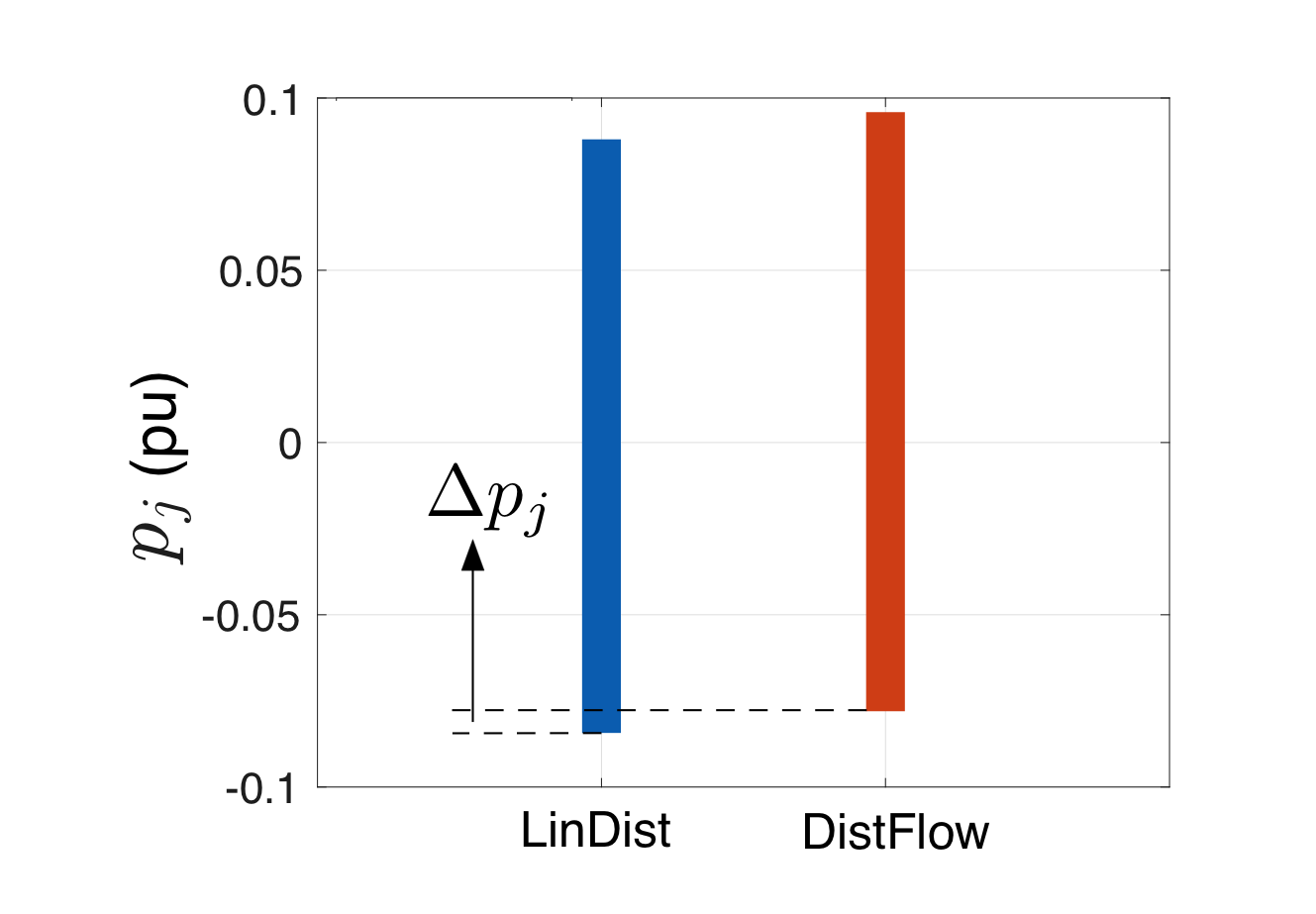}\label{fig:feas_region}}
\caption{(a) Comparison in change in node voltage with change in real power set-points between \textit{DistFlow} and \textit{LinDist} for two node system~\cite{nazir2019convex}. The two models result in different voltages and if design is based on linear model then network voltages limit will be violated if $p_j$ is operated at its lower limit, (b) Comparison of feasible region obtained from \textit{LinDist} and \textit{DistFlow} models for two node model. }
\end{figure}

The \textit{LinDist} model is obtained by neglecting the $l_{ij}$ term in~\eqref{eq:volt_rel}-\eqref{eq:curr_rel}. Consider a simple two-node system, where the active power injection, $p_j$, is varied at the load node and then observe the corresponding changes in voltage magnitude, $|V_j|$. As can be seen from Fig.~\ref{fig:QvsV}, the voltages obtained from the \textit{DistFlow} model (that represents the actual non-linear physics) and the \textit{LinDist} model (that sets $l_{ij}=0 \ \forall ij$) are identical at $p_j=0$. However, as we move $p_j$ away from this point, \textit{LinDist} becomes less accurate and this is especially the case towards the edge of the allowable power injections range (when voltages are near their upper and lower limits). Correspondingly, the models arrive at different conclusions when it comes to estimating the maximum and minimum active power injections as shown in Fig.~\ref{fig:feas_region}. Specifically, the \textit{LinDist} model overestimates the lower bound of the active power injection, which can in practice lead to voltage violations due to the model mismatch between \textit{LinDist} and \textit{DistFlow} models. Thus, linear models can clearly overestimate a network's capacity for flexibility and are hence not well-suited for grid-aware resource allocation problems.

\subsubsection{Shortcomings of convex relaxations}
\begin{figure}[!t]
\centering
\includegraphics[width=0.4\textwidth]{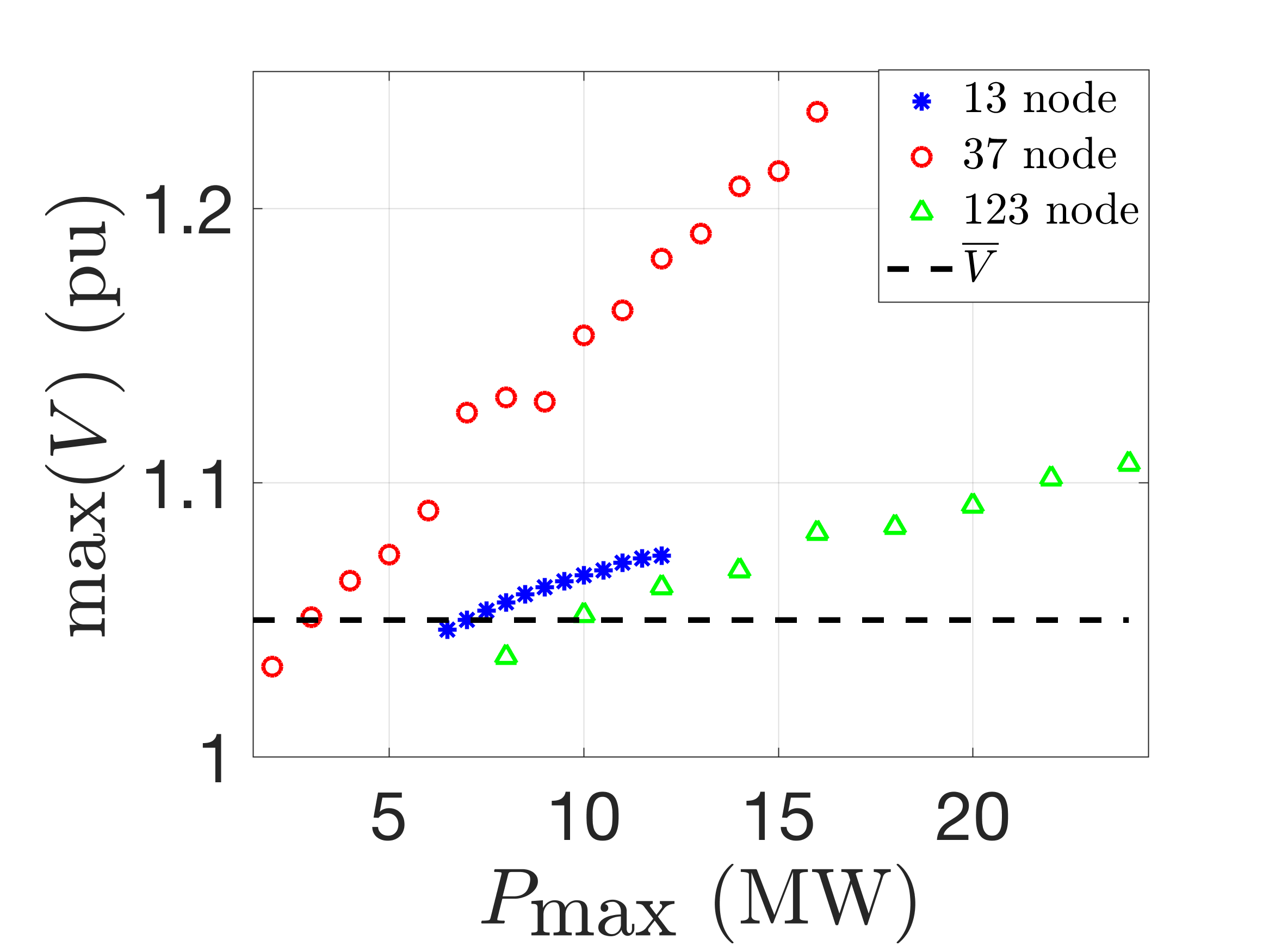}
\caption{\label{fig:HC_convrel} Plot of maximum nodal voltage with $P_{\text{max}}$ shows that the predicted hosting capacity obtained with a convex relaxation in (CR) is much higher than actually achievable. The figure shows the solutions obtained from the convex relaxation method (SOCP) are infeasible at the determined hosting capacity. As the maximum hosting capacity, $P_{\text{max}}$,  is reduced, a feasible solution that satisfies all network constraints is found.}.
\end{figure}

In this section, we expand on the above discussion to illustrate how approaches based on convex relaxations also can overestimate flexible capacity of a network, which can negatively impact reliability. Specifically, we consider the conventional SOCP relaxation of \textit{DistFlow}, obtained by replacing the nonlinear equality constraint in~\eqref{eq:curr_rel} with a convex inequality constraint. We also impose an upper bound on the maximum hosting capacity to determine the hosting capacity for which the AC network limits are satisfied. Solving the SOCP for three different distribution networks (13-node, 37-node and 123-node) results in Fig.~\ref{fig:HC_convrel}, which shows the maximum network voltage magnitude for different values of $P_{\text{max}}$. The examples show that the hosting capacity obtained from an SOCP solution will exceed the voltage limits unless $P_{\text{max}}$ is reduced significantly. For example, the SOCP formulation for the 37-node feeder will erroneously predict that the total allowable active power injections in the feeder exceeds 15MW. However, total allowable injections beyond 2.5MW can exceed the voltage limits at one or more nodes in the network. This is because a convex relaxation represents an \textit{outer approximation} of the original \textit{DistFlow} model and, hence, solutions obtained are not guaranteed to be physically realizable, i.e., they may violate network constraints despite being feasible in convex SOCP. 
For further discussion on convex OPF techniques, please see~\cite{gan2014exact,nick2017exact,nazir2020gridaware}. 

These shortcomings of the linear and the convex relaxation based models highlight the need for an inner approximation that provides admissibility guarantees. A convex form of an inner approximation would result in a scalable implementation that inherits the benefits of convex optimization with network admissibility guarantees. Recent works in literature have developed such a convex inner approximation of the power flow equations~\cite{nazir2020gridaware}, which can be used to determine a feeder's maximum, network-admissible active power nodal injections (i.e., nodal hosting capacity). A summary of the convex inner approximation formulation from~\cite{nazir2020gridaware} will be provided next.


\section{Convex inner approximations preliminaries}
In order to obtain a convex inner approximation (CIA) formulation first we define the vector form of~\eqref{eq:volt_rel}-\eqref{eq:reac_power_rel} as:

\begin{align}
    V=&v_{\text{0}}\mathbf{1}_N+M_{\text{p}}p+M_{\text{q}}q-Hl,\label{eq:final_volt_rel}\\
    P=&Cp-D_{\text{R}}l, \qquad Q=Cq-D_{\text{X}}l,\label{eq:P_relation}
\end{align}
where $P:=[P_{ij}]_{(i,j) \in \mathcal{L}} \in \mathbb{R}^N$, 
$Q:=[Q_{ij}]_{(i,j)\in \mathcal{L}}\in \mathbb{R}^N$, $V:=[v_i]_{i \in \mathcal{N}}\in \mathbb{R}^N$, $p:=[p_i]_{i \in \mathcal{N}}\in \mathbb{R}^N$, $p_{\text{g}}:=[p_{\text{g},i}]_{i \in \mathcal{N}}\in \mathbb{R}^N$, $P_{\text{L}}:=[P_{\text{L},i}]_{i \in \mathcal{N}} \in \mathbb{R}^N$, $q:=[q_i]_{i \in \mathcal{N}}\in \mathbb{R}^N$, $Q_{\text{L}}:=[Q_{\text{L},i}]_{i \in \mathcal{N}} \in \mathbb{R}^N$, and  $l:=[l_{ij}]_{(i,j) \in \mathcal{L}}\in \mathbb{R}^N$
and matrices $R:=\text{diag}\{r_{ij}\}_{(i,j)\in \mathcal{L}}\in \mathbb{R}^{N\times N}$, $X:=\text{diag}\{x_{ij}\}_{(i,j)\in \mathcal{L}}\in \mathbb{R}^{N\times N}$, $Z^2:=\text{diag}\{z_{ij}^2\}_{(i,j)\in \mathcal{L}}\in \mathbb{R}^{N\times N}$, and $A:=[0_N \quad I_N]B-I_N$ and $B\in \mathbb{R}^{(N+1)\times N}$ is the \textit{incidence matrix} of $\mathcal{G}$ relating the branches in $\mathcal{L}$ to the nodes in $\mathcal{N}\cup \{0\}$, where $I_N$ is the $N\times N$ identity matrix and $0_N$ is a column vector of $N$ rows. Also matrices $M_{\text{p}}:=2C^TRC$, $\quad$ $M_{\text{q}}:=2C^TXC$, $\quad$ $H:=C^T(2(RD_{\text{R}}+XD_{\text{X}})+Z^2)$ and $C:=(I_N-A)^{-1}$, $D_{\text{R}}:=(I_N-A)^{-1}AR$, and $D_{\text{X}}:=(I_N-A)^{-1}AX$ describe the network topology and impedance parameters.
Further details about this formulation of the \textit{DistFlow} equations can be found in~\cite{heidari2017non,nazir2019voltage,nazir2020gridaware}.

Based on the linear model of the \textit{DistFlow} equations in~\eqref{eq:final_volt_rel}-\eqref{eq:P_relation} and considering $l_{\text{lb}}$ and $l_{\text{ub}}$ as the new proxy variables for lower and upper bounds on $l$. Then, we can define the corresponding upper $(.)^+$ and lower $(.)^-$ bounds of $P$, $Q$ and $V$ as follows:

\begin{subequations}\label{eq:CIA_bounds}
\begin{align}
    P^+:=&Cp-D_{\text{R}}l_{\text{lb}} \label{eq:P_relation_1}\\
    P^-:=&Cp-D_{\text{R}}l_{\text{ub}} \label{eq:P_relation_2}\\
    Q^+:=&Cq-D_{\text{X}_+}l_{\text{lb}}-D_{\text{X}_-}l_{\text{ub}}\label{eq:Q_relation_1}\\
    Q^-:=&Cq-D_{\text{X}_+}l_{\text{ub}}-D_{\text{X}_-}l_{\text{lb}}\label{eq:Q_relation_2}\\
    V^+:=&v_{\text{0}}\mathbf{1}_n+M_{\text{p}}p+M_{\text{q}}q-H_+l_{\text{lb}}-H_-l_{\text{ub}}\label{eq:V_relation_1}\\
    V^-:=&v_{\text{0}}\mathbf{1}_n+M_{\text{p}}p+M_{\text{q}}q-H_+l_{\text{ub}}-H_{-} l_{\text{lb}}, \label{eq:V_relation_2}
\end{align}
\end{subequations}
where $D_{\text{X}_+}$ and $H_+$ include the non-negative elements of $D_{\text{X}}$ and $H$, respectively, and $D_{\text{X}_-}$ and $H_-$ are the corresponding negative elements. These bounds $l_\text{lb}, l_\text{ub}$ allow us to neglect the non-linear term in~\eqref{eq:curr_rel}. Next we will obtain the expression for these bounds on $l$.

Based on any nominal operating point $x^0_{ij}:=\text{col}\{P_{ij}^0, Q_{ij}^0, v_j^0\} \in\mathbb{R}^3$, the second-order approximation for~\eqref{eq:curr_rel} can be expressed as
\begin{align}\label{eq:T_exp}
    l_{ij} & \approx l_{ij}^0 + \mathbf{J}_{ij}^\top \mathbf{\delta}_{ij} +\frac{1}{2}\mathbf{\delta}_{ij}^\top \mathbf{H}_{\text{e},ij} \mathbf{\delta}_{ij},
\end{align}
where $l_{ij}^0 := l_{ij}(x_{ij}^0)$ are squared branch currents, $\mathbf{\delta}_{ij}:=[P_{ij}-P_{ij}^0, \, Q_{ij}-Q_{ij}^0, \, v_j-v_j^0]^\top$ and Jacobian $\mathbf{J}_{ij}\in \mathbb{R}^{3\times 1}$ and Hessian $\mathbf{H}_{\text{e},ij}\in \mathbb{R}^{3\times 3}$ are defined for each branch $(i,j)$ as detailed in~\cite{nazir2020gridaware}.
 Furthermore, since $\mathbf{H}_{\text{e},ij}$ is provably positive semi-definite per~\cite{nazir2019voltage} the second-order approximation in~\eqref{eq:T_exp} can be used to derive lower and upper bounds of $l_{ij}$ for all $(i,j)\in \mathcal{L}$. This is described as follows:
 \begin{align}
    l_{ij}=|l_{ij}| & \approx |l_{ij}^0 + \mathbf{J}_{ij}^\top\mathbf{\delta}_{ij}+\frac{1}{2}\mathbf{\delta}_{ij}^\top \mathbf{H}_{\text{e},ij} \mathbf{\delta}_{ij}| \\
    & \le |l_{ij}^0| + |\mathbf{J}_{ij}^\top\mathbf{\delta}_{ij}|+|\frac{1}{2}\mathbf{\delta}_{ij}^\top \mathbf{H}_{\text{e},ij} \mathbf{\delta}_{ij}| \\
    & \le l_{ij}^0 + \max\{2|\mathbf{J}_{ij}^\top\mathbf{\delta}_{ij}|,|\mathbf{\delta}_{ij}^\top \mathbf{H}_{\text{e},ij} \mathbf{\delta}_{ij}|\}\label{eq:l_upper_1}\\
    \Rightarrow  l_{ij} &
     \le l_{ij}^0 + \max\{2|\mathbf{J}_{ij+}^\top\mathbf{\delta}_{ij}^++\mathbf{J}_{ij-}^\top\mathbf{\delta}_{ij}^-|,\mathbf{\psi}_{ij}\} \le  l_{\text{ub},ij}, \label{eq:l_upper}\\
     \text{and } l_{ij} & \ge l_{ij}^0 + \mathbf{J}_{ij+}^\top\mathbf{\delta}_{ij}^-+\mathbf{J}_{ij-}^\top\mathbf{\delta}_{ij}^+ =: l_{\text{lb},ij}, \label{eq:l_lower}
\end{align}
where $\mathbf{J}_{ij+}$ and $\mathbf{J}_{ij-}$ includes the positive and negative elements of $\mathbf{J}_{ij}$, $\mathbf{\delta}_{ij}^+ := \mathbf{\delta}_{ij}(P_{ij}^+, Q_{ij}^+, v_j^+,x_{ij}^0)$ and $\mathbf{\delta}_{ij}^-:=\mathbf{\delta}_{ij}(P_{ij}^-, Q_{ij}^-, v_j^-,x_{ij}^0)$, and $\mathbf{\psi}_{ij}:=\max\{(\mathbf{\delta}_{ij}^{+,-})^\top \mathbf{H}_{\text{e},ij} (\mathbf{\delta}_{ij}^{+,-})\}$, which represents the largest of eight possible combinations of $P/Q/v$ terms in $\mathbf{\delta}_{ij}$ with mixed $+,-$ superscripts. A detailed description of these terms together with an analysis on the accuracy of this second-order Taylor series approximation can be found in~\cite{nazir2019voltage,nazir2020gridaware}. 

This inner approximation of the \textit{DistFlow} equations~\eqref{eq:P_relation_1}-\eqref{eq:V_relation_2},\eqref{eq:l_upper},\eqref{eq:l_lower} is provably convex and will henceforth be referred to as the convex inner approximation (CIA). This CIA can be employed to obtain a hyper-rectangle characterization of the set of active power nodal injections  where the corner points represent the feeders nodal capacity limits for the flexibility. Further analysis in~\cite{nazir2020gridaware} also proves the admissibility of the range  $\Delta p_{\text{g}} := [p_{\text{g}}^-, p_{\text{g}}^+]$ obtained through the CIA formulation and also presents an iterative algorithm that increases the size of the admissible set. The next section will lay out the market clearing problem to allocate flexibility amongst Aggregators based on the CIA formulation.



\section{Market-based grid-aware formulation}
\begin{figure}[!t]
\centering
\includegraphics[width=0.4\textwidth]{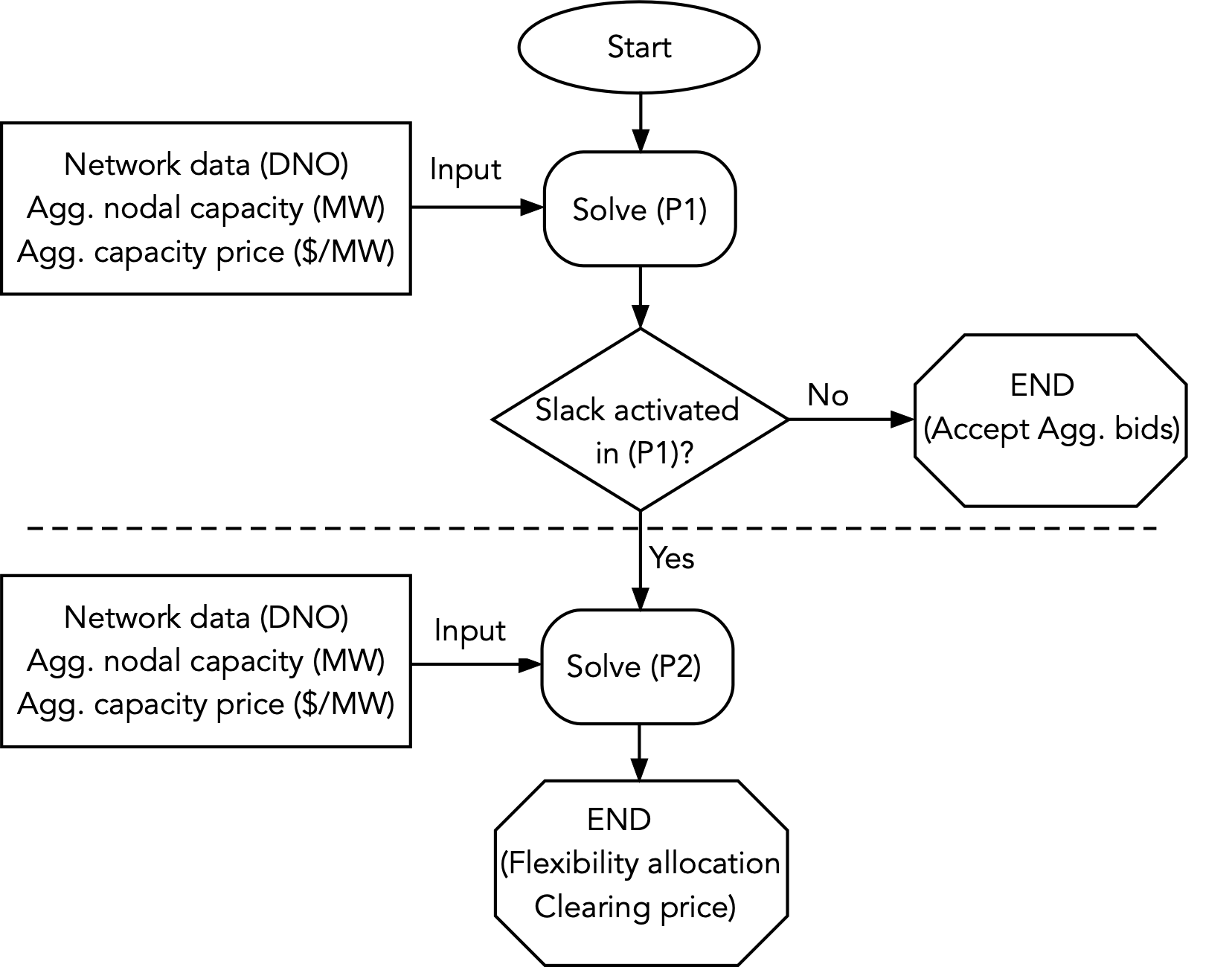}
\caption{\label{fig:market_flow} Illustration of the DSO's market clearing mechanism for grid-aware allocation of Aggregator flexibility.}
\end{figure}

To formulate the market-clearing problem, the DSO receives capacity and price bids from the Aggregators. The capacity bids represent an upper bound on the Aggregator's available flexibility at each node and the price bid represents an upper bound on acceptable capacity charges at each nodes. Based on these bids, the DSO allocates flexible capacity to Aggregators at each node in the distribution network. The DSO's market-clearing problem  is implemented as a two-step market optimization problem. The flowchart in Fig.~\ref{fig:market_flow} illustrates the two-step approach. In Step~1, the DSO checks feasibility of hosting all Aggregators' flexible capacity. If congestion is expected to result from any Aggregator dispatch, Step~2 is used to prioritize allocated flexibility based on the Aggregators' price bid, which represent their maximum access charges. This two-step approach is detailed next.

\subsection{Step~1: Network feasibility certificate}
In Step~1, the DSO first checks the feasibility of network to host the total Aggregators' flexibility (MW) and does not require prioritizing any node or Aggregator and, thus, does not require the price element from the Aggregators' bids. To determine the impact of active power injections, the DSO leverages the CIA presented above to formulate an inner approximation of the non-convex (NLP) problem. This is illustrated below in (P1). That is, in Step~1, the DSO solves (P1) and, if no slack variables are activated (i.e., $s_i^*=0 \,\,\forall i$), the DSO informs the Aggregators that all flexible capacity can be accommodated by the network and no premium access charge are transferred to the Aggregators.
This optimization problem can be formulated as:
{\color{black}
\begin{subequations}
\begin{align}
\text{(P1)}  \quad s_i^*=\arg\min_{s_i^+} \  \sum_{i=1}^{N}s_i^+&\\
\text{s.t.}  \quad  \eqref{eq:P_relation_1}-\eqref{eq:V_relation_2},\eqref{eq:l_upper},\eqref{eq:l_lower}\\
    p_i=p_{\text{g},i}^+-P_{\text{L},i}-s_i^+ & \quad \forall i\in \mathcal{N} \\
    \underline{v_i}\le v_i \le \overline{v_i} & \quad \forall i\in \mathcal{N}   \label{eq:P1_V}\\
     l_{\text{ub},ij} \le \overline{l_{ij}} & \quad \forall (i,j)\in \mathcal{L} \label{eq:P1_lmax}\\
    s_i^+ \ge 0 & \quad \forall i\in \mathcal{N}
\end{align}
\end{subequations}
}
where $p_{\text{g},i}^+:=\sum_mp^{\text{bid}}_{m,i}$, with $p^{\text{bid}}_{m,i}$ being the flexible capacity bid (MW) of the $mth$ Aggregator at node $i$, where $m\in \mathcal{M}$, with $\mathcal{M}$ being the set of Aggregators. Also, let $k_{m,i}$ be the price bid (\$/MW) provided by the $m$-th Aggregator to the DSO to define its upper bound on acceptable network access charges at node $i$.
Now, if (P1) predicts that any slack variables will become active (i.e., any $s_i^*>0$), then some Aggregator flexibility must be curtailed and, thus, we require a grid-aware allocation that prioritizes nodal flexibility in certain locations and from certain Aggregators (who are willing to pay higher access fees). Thus, we need to leverage the available pricing information from each Aggregator. This prioritization is described in Step~2 in the following section.

\subsection{Step~2: Network-aware optimal allocation}
If the DSO predicts network congestion by solving (P1), then the DSO needs to allocate flexible capacity optimally based on price information. This is the goal of Step~2, where we solve for an optimal capacity allocation of the Aggregators based on both the capacity (MW) and price  bids (\$/MW). Thus, Step~2 determines the maximum DNO access charge that an Aggregator is willing to pay.
This represents a generalization of the author's prior work on  CIA~\cite{nazir2020gridaware}, where the formulation now support multiple Aggregators competing for nodal hosting capacity via flexibility capacity bids (MW and \$/MW). Specifically, the Aggregator capacities feature in the energy balance constraints while the price bids enter in the objective function to prioritize aggregator flexibility across the network.

 Based on the price bids and flexibility offered by various Agregators, the DSO solves a centralized market optimization problem that determines the admissible flexibility range and the market clearing price for each Aggregator at each node, i.e., $\Delta p_{m,i}=p^+_{m,i}-p^-_{m,i}$. The optimization problem is formulated as:
 \begin{subequations}
\begin{align}
   \text{(P2)}  \quad p^+_{m,i}=\arg\max_{ p_{m,i}}\sum_m\sum_i k_{m,i} p_{m,i}\\
    s.t. \ p^{\text{bid}}_{m,i}- p_{m,i}\ge 0 \quad \forall m\in \mathcal{M}, \forall i\in \mathcal{N} \\
    p_i=\sum_m p_{m,i}-P_{\text{L},i} \quad \forall i\in \mathcal{N}\\
    \eqref{eq:P_relation_1}-\eqref{eq:V_relation_2},\eqref{eq:l_upper},\eqref{eq:l_lower}\\
    \underline{V}\le     V^-(p,q) \quad V^+(p,q)\le &\overline{V} \\
     l_{\text{ub}}\le  \overline{l}.
\end{align}
\end{subequations}
\textcolor{black}{ A similar optimization problem needs to be solved for the lower ranges, $p_{\text{m},i}^-$.}
To illustrate this two-step market mechanism, next we present case studies on the modified IEEE-37 node distribution network with multiple Aggregators.

\subsubsection{Case study~1: Feeder-level prices per Aggregator}
\begin{figure}
    \centering
    \subfloat{
    \includegraphics[width=0.48\linewidth]{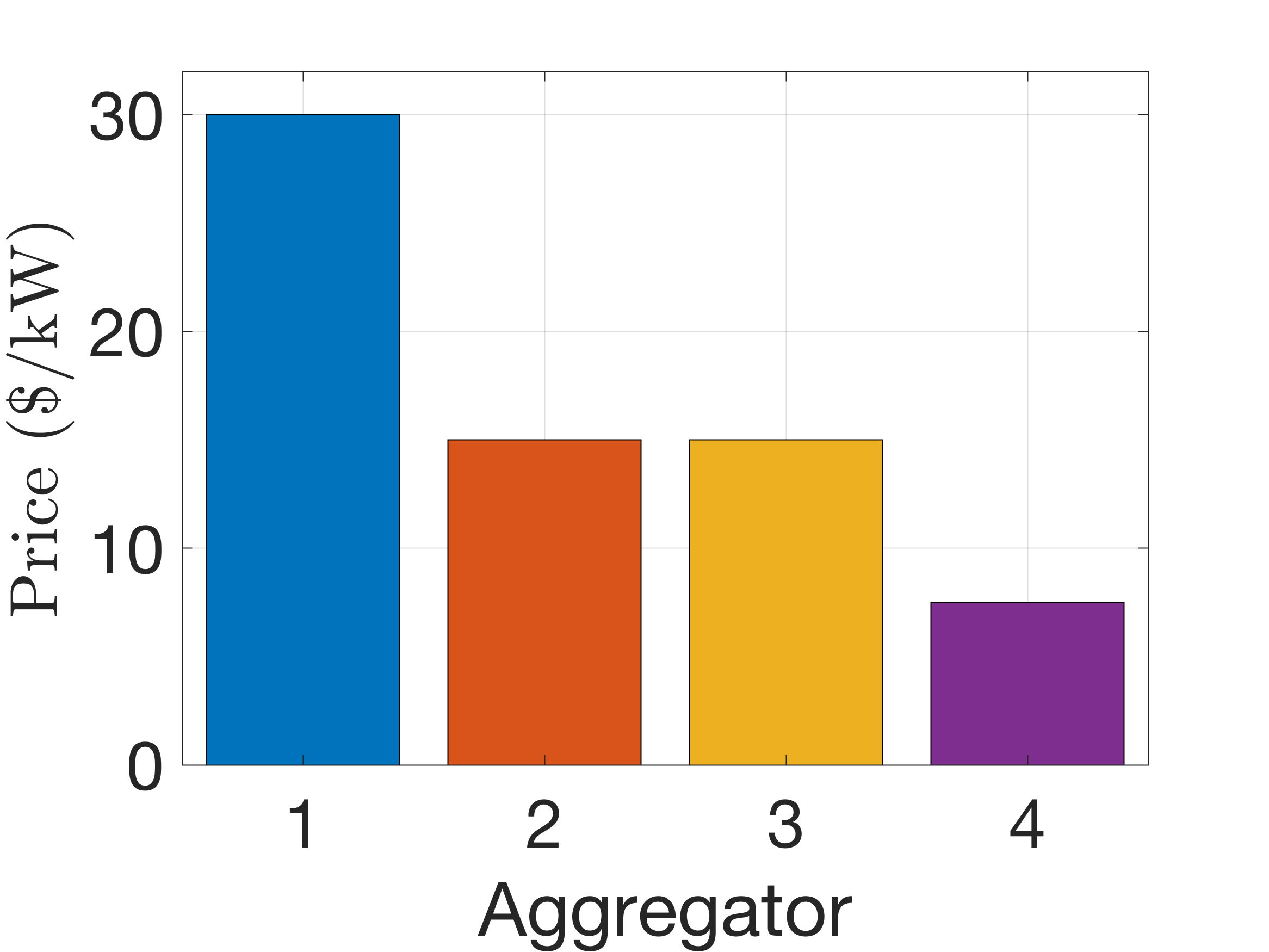}\label{fig:Agg_price}}
    \subfloat{
    \includegraphics[width=0.48\linewidth]{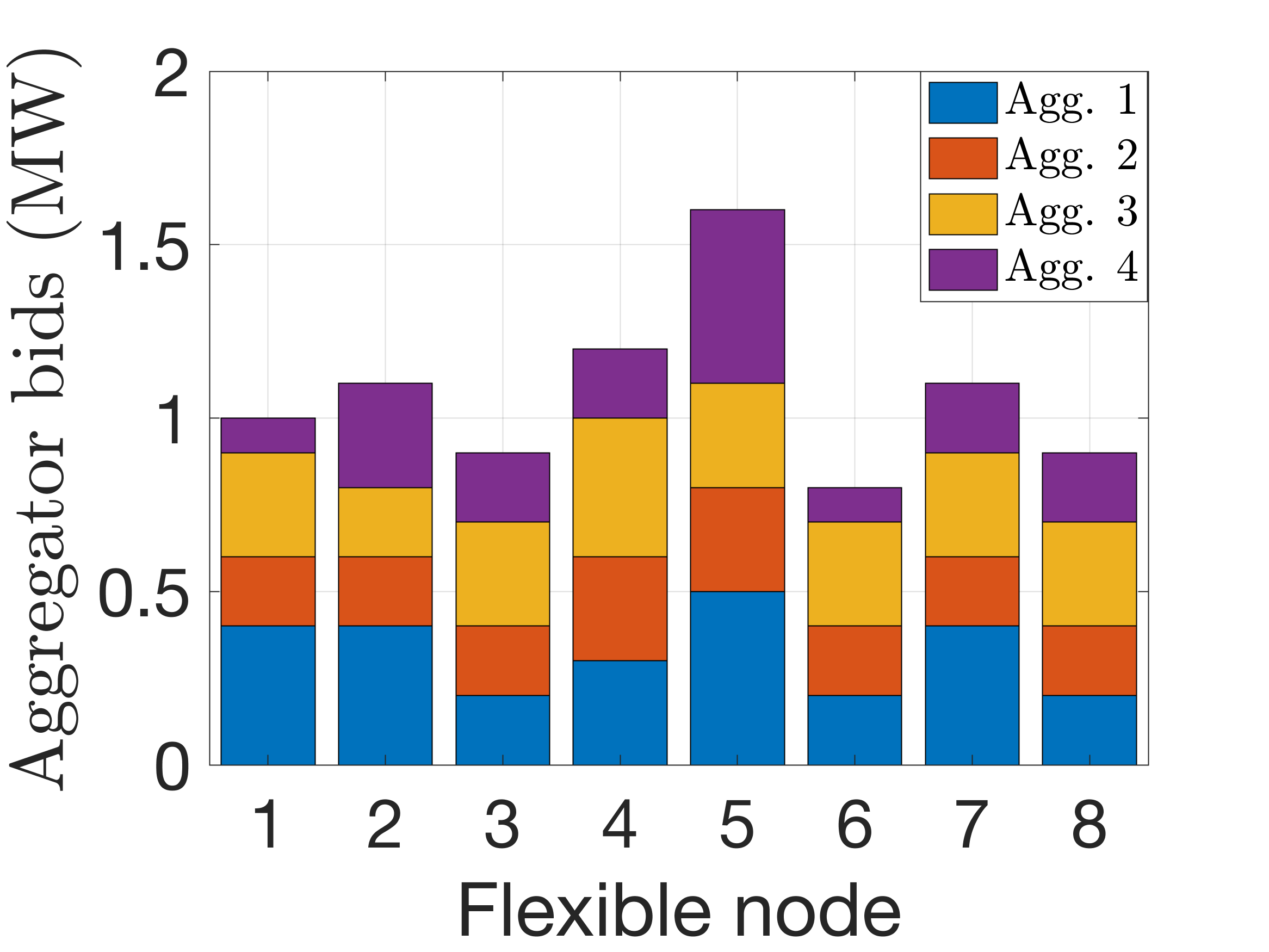}\label{fig:Agg_bids}}
\caption{(a) Day-ahead Aggregator prices in a distribution feeder (assuming Aggregators have same prices at different nodes), (b) Day-ahead Aggregator flexibility bids into the the DSO clearing market. }
\end{figure}

To showcase the market-clearing problem (P1) and (P2), we utilize a modified IEEE-37 node system from~\cite{baker_cognizant} and consider a scenario where each Aggregator's bid  consists of nodal capacities (MW) and a fixed capacity price ($\$/MW$) for all nodes as shown in Fig.~\ref{fig:Agg_price}. The price bid values are adopted from commercial demand charges based on available cost-optimization models~\cite{mclaren2017identifying}. \textcolor{black}{In this case study we only consider the upper range of flexibility as depicted in Fig.~\ref{fig:Agg_bids}.} Based on the Aggregator bids, the DSO first solves (P1) and checks if $||s_i^*||_{\infty} > \epsilon $, where tolerance $\epsilon$ chosen to be 10 W. In this case study, (P1) yields $||s_i^*||_{\infty}=0.9$ MW, which indicates that the Aggregators' flexibility could cause network congestion. Hence, the DSO proceeds to Step~2 and solves (P2), which uses the capacity prices to prioritize Aggregators access and satisfy network constraints.

\begin{figure}[!t]
    \centering
    \subfloat{
    \includegraphics[width=0.48\linewidth]{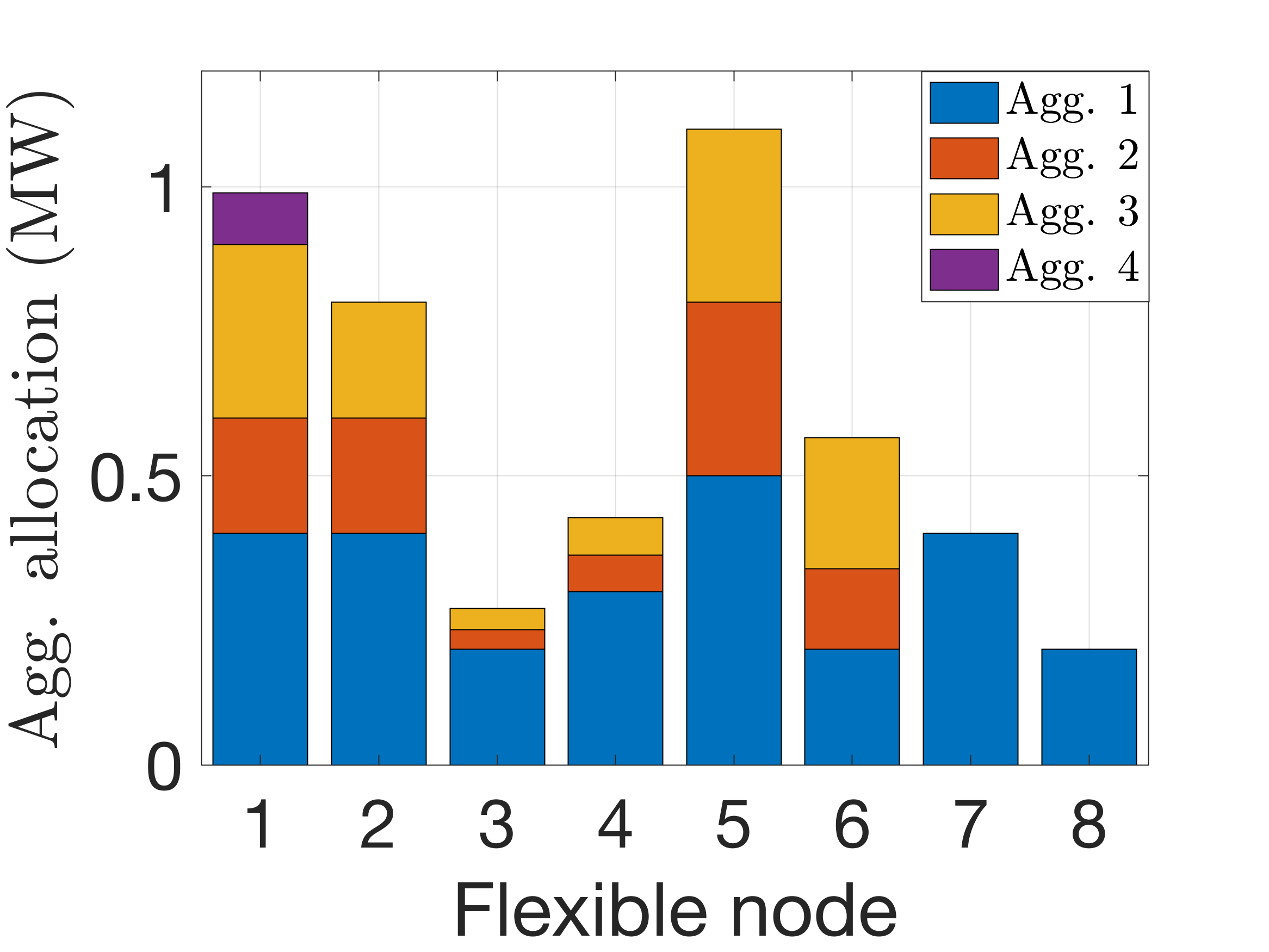}\label{fig:Agg_alloted}}
    \subfloat{
    \includegraphics[width=0.48\linewidth]{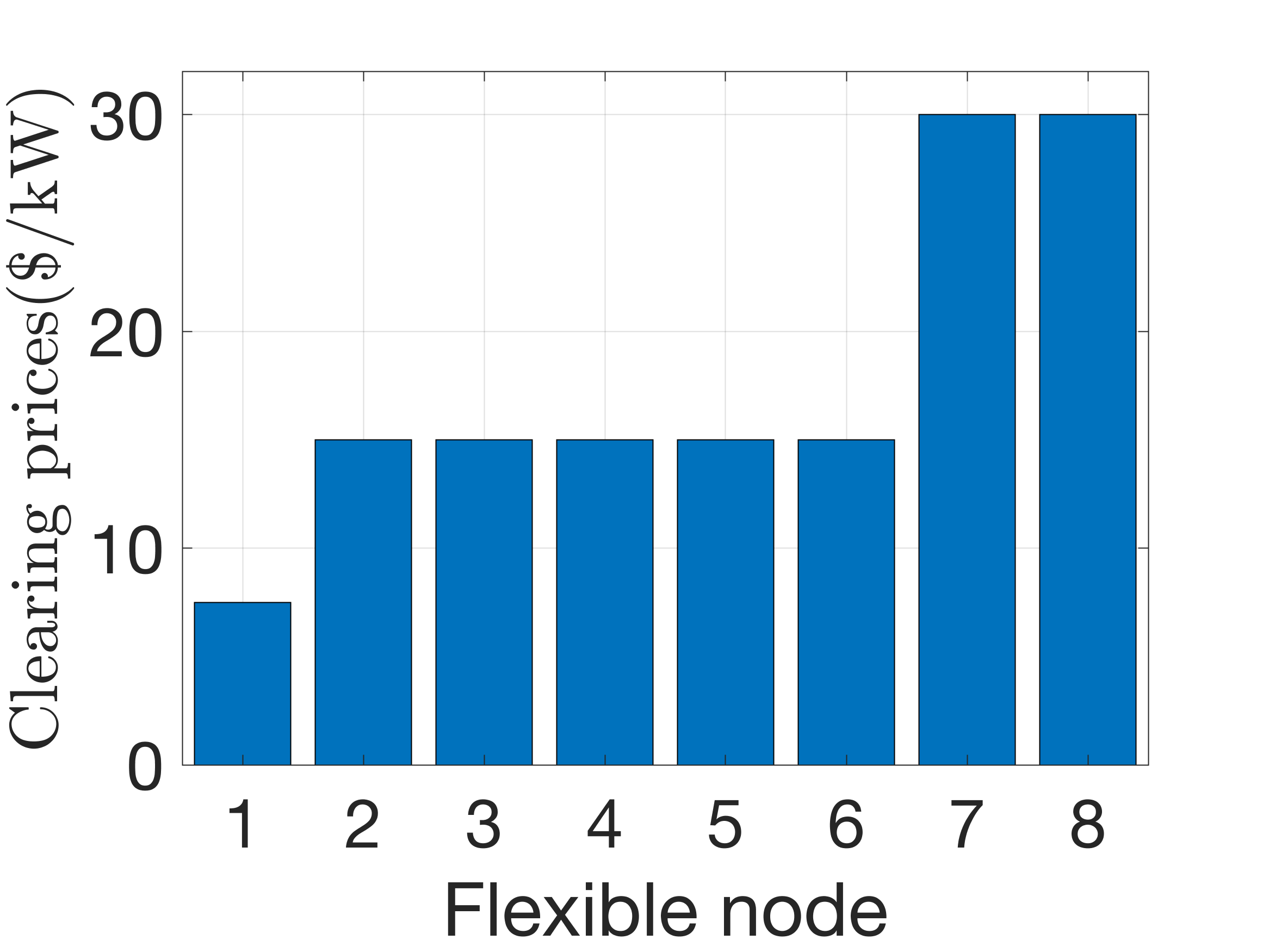}\label{fig:clearing_base}}
\caption{Case Study~1 results. (a)  Allocation of Aggregator flexibility after market clearing, (b) Clearing price after market allocation based on cheapest allocation at each node. }
\end{figure}

 Based on predicted congestion, DSO solves (P2) and the resulting flexibility allocation is shown in Fig.~\ref{fig:Agg_alloted} and highlights that Aggregators willing to pay less for hosting capacity are allocated less flexibility. Here, Aggregator~1 is prioritized due to its willingness to pay a higher price to access the network's capacity. This point is further illustrated in Table~\ref{table_base} which shows what fraction of Aggregator and nodal flexibility is allocated by the DSO. The results show that nodes close to the substation get priority in capacity allocation.

After the Aggregators' flexible capacities are allocated by the DSO, we utilize a \textcolor{black}{price-clearing scheme where the lowest price at each node decides the market clearing price~\cite{sinha2008setting} as depicted} in Fig.~\ref{fig:clearing_base}. The figure shows three different clearing prices, with the nodes closer to the substation having a lower price because of less congestion at those nodes. The revenue earned by the DNO in this case is calculated according to~\eqref{eq:revenue} and comes out to $\$72,870$. 
\begin{align}\label{eq:revenue}
   \text{Total DNO Revenue}=\sum_i \sum_m k_{\text{c},i}p_{m,i}
\end{align}
where $k_{\text{c},i}$ is the clearing price at node $i$.

This case study employed feeder-level prices per Aggregator, which means that all nodes are priced the same by each Aggregator, which could represent some average customer or flexible load acquisition costs. In the next case study, we explore nodal-level capacity prices for the Aggregators, which is similar to the concept of distribution locational marginal prices (DLMPs)~\cite{Bai_DLMP}. However, unlike much of the literature on DLMPs, this work focuses on pricing access to network capacity (similar to demand charges) rather than incentivizing DER participation by broadcasting (variable) prices in the distribution network.

\begin{table}[!t]
\caption{Fraction of flexible capacity allocated across nodes and aggregators}
\label{table_base}
\centering
\begin{tabular}{rllllllll}
\hline
{Item} & {1} & {2} & {3} & {4} & {5} & {6} & {7} & {8} \\
\hline
Nodes & 0.99 & 0.73 & 0.3 & 0.36 & 0.69 & 0.71 & 0.36 & 0.22 \\
Agg. & 1 & 0.52 & 0.47 & 0.05 \\
\hline
\end{tabular}
\end{table}


\subsubsection{Case study 2: Nodal Aggregator pricing}

\begin{figure}[!t]
    \centering
    \subfloat{
    \includegraphics[width=0.48\linewidth]{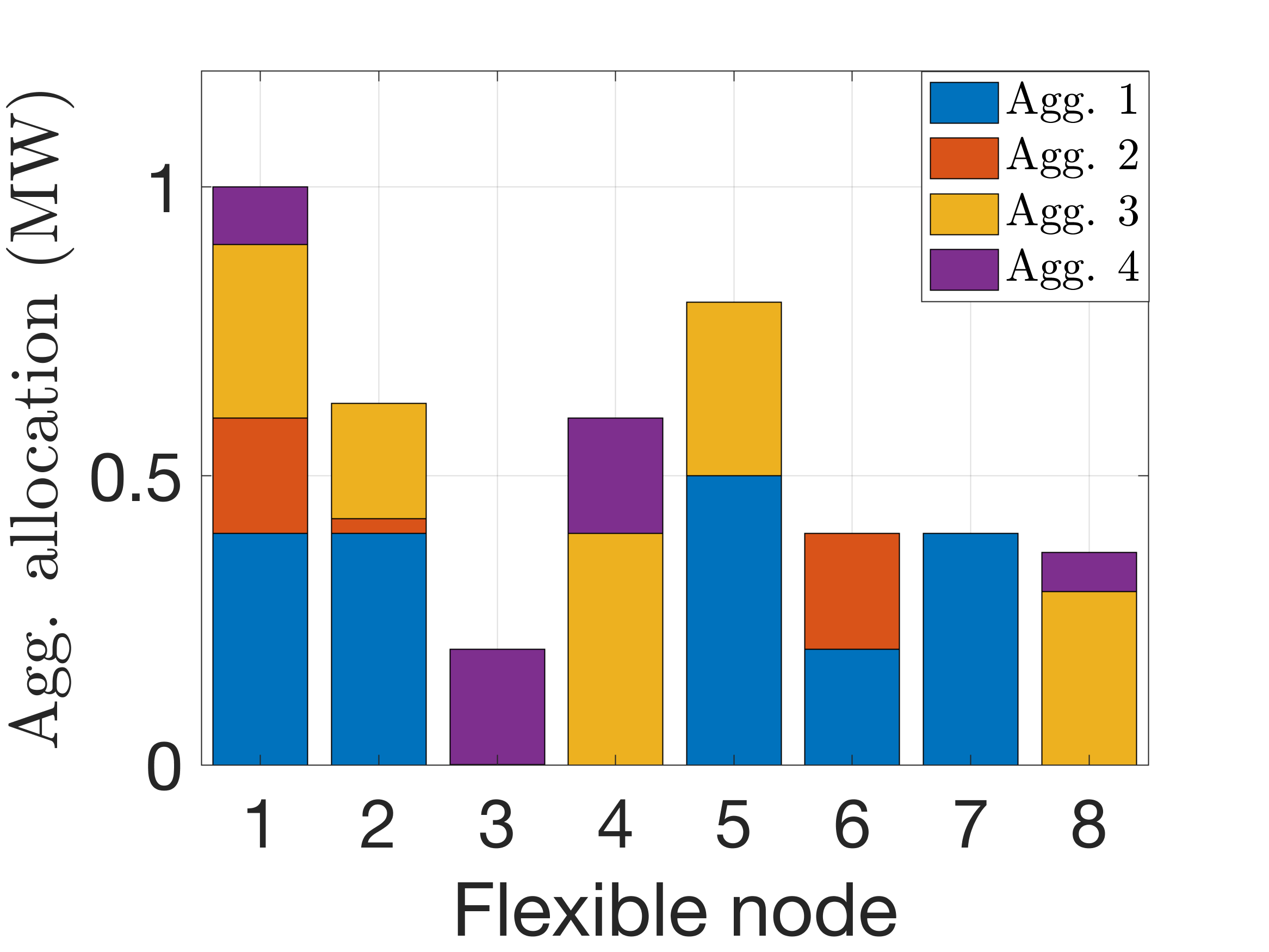}\label{fig:Agg_alloted_VP}}
    \subfloat{
    \includegraphics[width=0.48\linewidth]{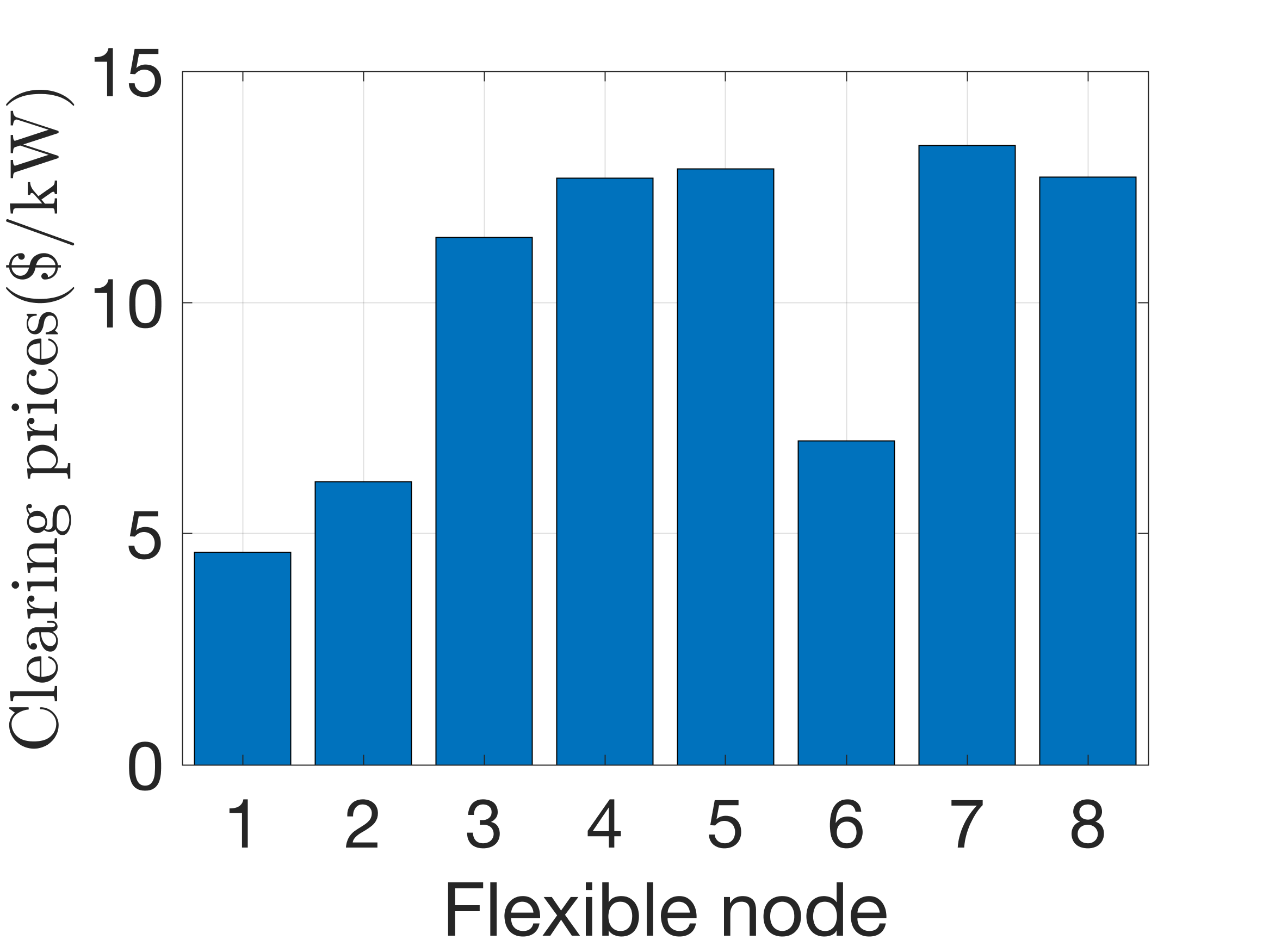}\label{fig:clearing_VP}}
\caption{Case Study~2 results. (a) Allocation of Aggregator flexibility after market clearing under nodal Aggregator pricing, (b) Clearing price after market allocation based on cheapest allocation at each node}
\end{figure}

In this case study, we consider a nodal-level Aggregator pricing structure as shown in Table~\ref{table_var} for the four Aggregators considered. \textcolor{black}{The price bid values are chosen randomly but are again based on the concept of commercial demand charges based on available cost-optimization models~\cite{mclaren2017identifying}.} Based on these prices, we solve the (P2) market optimization problem to determine the capacity allocations which are shown in Fig.~\ref{fig:Agg_alloted_VP}. Based on this allocation, we again compute the clearing price based on the least-price allocation of an Aggregator at each node, which is shown in Fig.~\ref{fig:clearing_VP}. The DNO revenue in this scenario is calculated to be $\$41,460$.

\begin{table}[!t]
\caption{Data for nodal-level aggregator capacity prices (\$/MW)}
\label{table_var}
\centering
\begin{tabular}{rllllllll}
\hline
{Node} & {1} & {2} & {3} & {4} & {5} & {6} & {7} & {8} \\
\hline
Agg. 1 & 9.8  & 10.9 & 1.2 & 2.5 & 12.9 & 13.5 & 13.4  & 4.1 \\
Agg. 2 & 4.6 & 6.1 & 5.6 & 1.6 & 3.60 & 7.0 & 4.6 & 2.6 \\
Agg. 3 & 29.1 & 9.0 & 6.8 & 14.0 & 21.0 & 2.7 & 0.9 & 15.6 \\
Agg. 4 & 13.2 & 1.6 & 11.4 & 12.7 & 1.3 & 4.8 & 2.5 & 12.7 \\
\hline
\end{tabular}
\end{table}

\section{Robust allocation of flexibility}
The analysis provided up to this point has been based on the assumption that the background demand is known in advance. However, in practice that may not be the case and it becomes important for the DSO to consider the uncertainties in forecasting background demand. To address this problem, we augment optimization formulations (P1) and (P2).

\subsection{Robust formulation}
In this method, we consider the uncertainty in background demand to be bounded and obtain a formulation that is robust to the worst-case of the bounded uncertainty. The robust versions of (P1) is formulated as shown below:
{\color{black}
\begin{subequations}
\begin{align}
\text{(ROB-P1)}  \quad s_i^*=\arg\min_{s_i^+} \  \sum_{i=1}^{N}s_i^+&\\
\text{s.t.}  \quad  \eqref{eq:P_relation_1}-\eqref{eq:V_relation_2},\eqref{eq:l_upper},\eqref{eq:l_lower}\\
    p_i=p_{\text{g},i}^+-P_{\text{L},i}-s_i^++D_i^+ & \quad \forall i\in \mathcal{N}     \label{eq:P1_disturbance} \\
    \underline{v_i}\le v_i \le \overline{v_i} & \quad \forall i\in \mathcal{N}   \label{eq:RP1_V}\\
     l_{\text{ub},ij} \le \overline{l_{ij}} & \quad \forall (i,j)\in \mathcal{L} \label{eq:RP1_lmax}\\
    s_i^+ \ge 0 & \quad \forall i\in \mathcal{N}
\end{align}
\end{subequations}
}
\textcolor{black}{where $D_i^+$ in~\eqref{eq:P1_disturbance} is the upper robust bound on the uncertainty in background demand at node $n$, which is considered alongside the upper flexibility range $p_{\text{g},i}^+$, respectively.} A similar problem would consider the lower robust bound on uncertainty $D_i^-$ in order to determine the feasibility of the lower flexibility range $p_{\text{g},i}^-$. In a similar fashion, the robust (P2) is formulated as follows:
\begin{subequations}
\begin{align}
    \text{(ROB-P2)}  \quad p^+_{m,i}=\arg\max_{ p_{m,i}}\sum_m\sum_i k_{m,i} p_{m,i}\\
    s.t. \ p^{\text{bid}}_{m,i}- p_{m,i}\ge 0 \quad \forall m\in \mathcal{M}, \forall i\in \mathcal{N} \\
    p_i=\sum_m p_{m,i}-P_{\text{L},i}+D_i^+(D_i^-) \quad \forall i \in \mathcal{N}\label{eq:P2_disturbance}\\
    \eqref{eq:P_relation_1}-\eqref{eq:V_relation_2},\eqref{eq:l_upper},\eqref{eq:l_lower}\\
     \underline{V}\le     V^-(p,q) \quad V^+(p,q)\le &\overline{V} \\
     l_{\text{ub}}\le  \overline{l}
\end{align}
\end{subequations}

 \textcolor{black}{A similar optimization problem would consider the lower robust bound on uncertainty $D_i^-$ in order to determine the lower flexibility allocation $p_{\text{m},i}^-$. This formulation ensures that the market allocation $p_{\text{m},i}^+(p_{\text{m},i}^-)$ accounts for the uncertainty $D_i^+(D_i^-)$.}

\subsubsection{Simulations under uncertainty}
\begin{figure}[!t]
    \centering
    \subfloat{
    \includegraphics[width=0.48\linewidth]{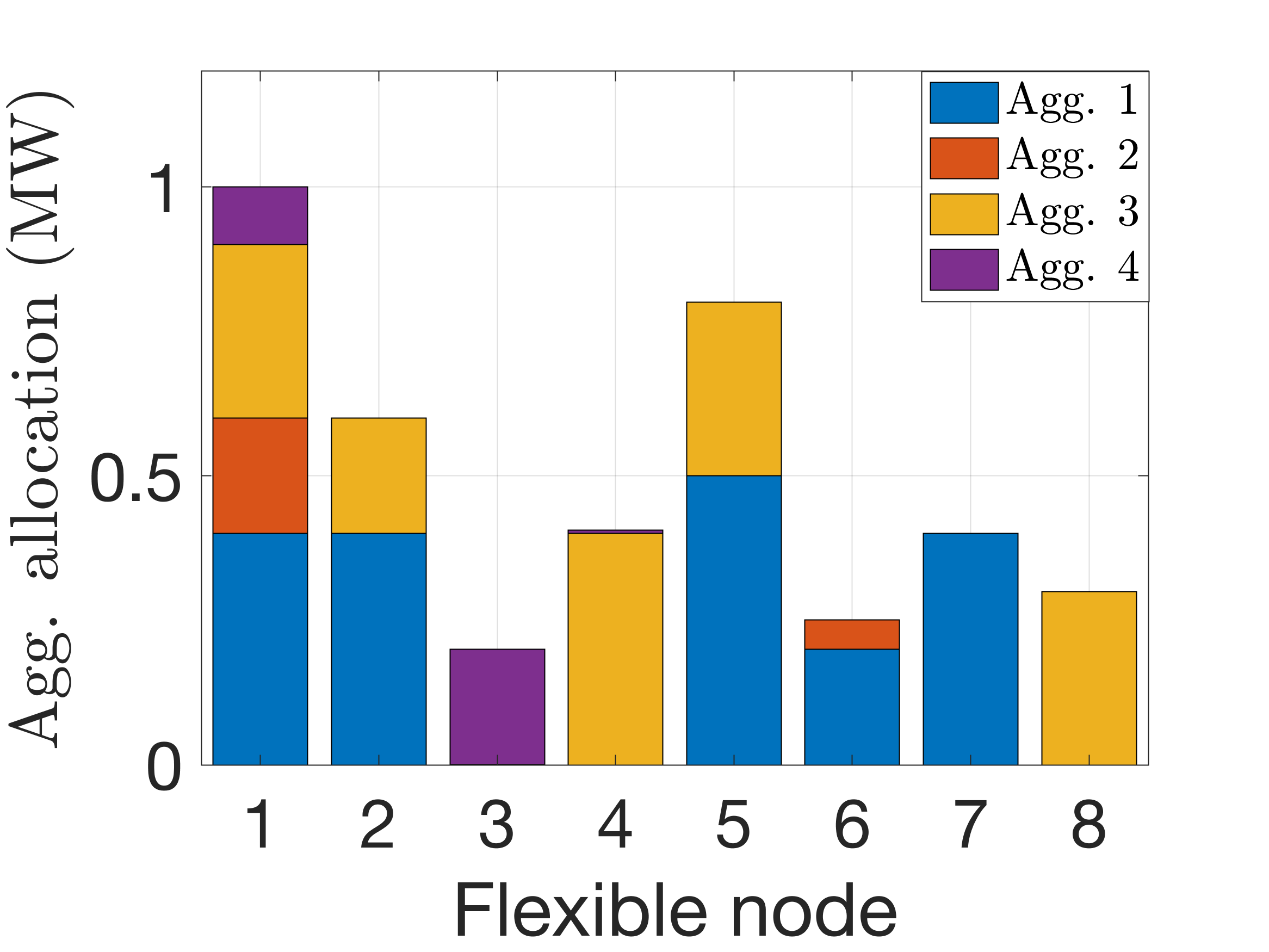}\label{fig:Agg_alloted_uncer}}
    \subfloat{
    \includegraphics[width=0.48\linewidth]{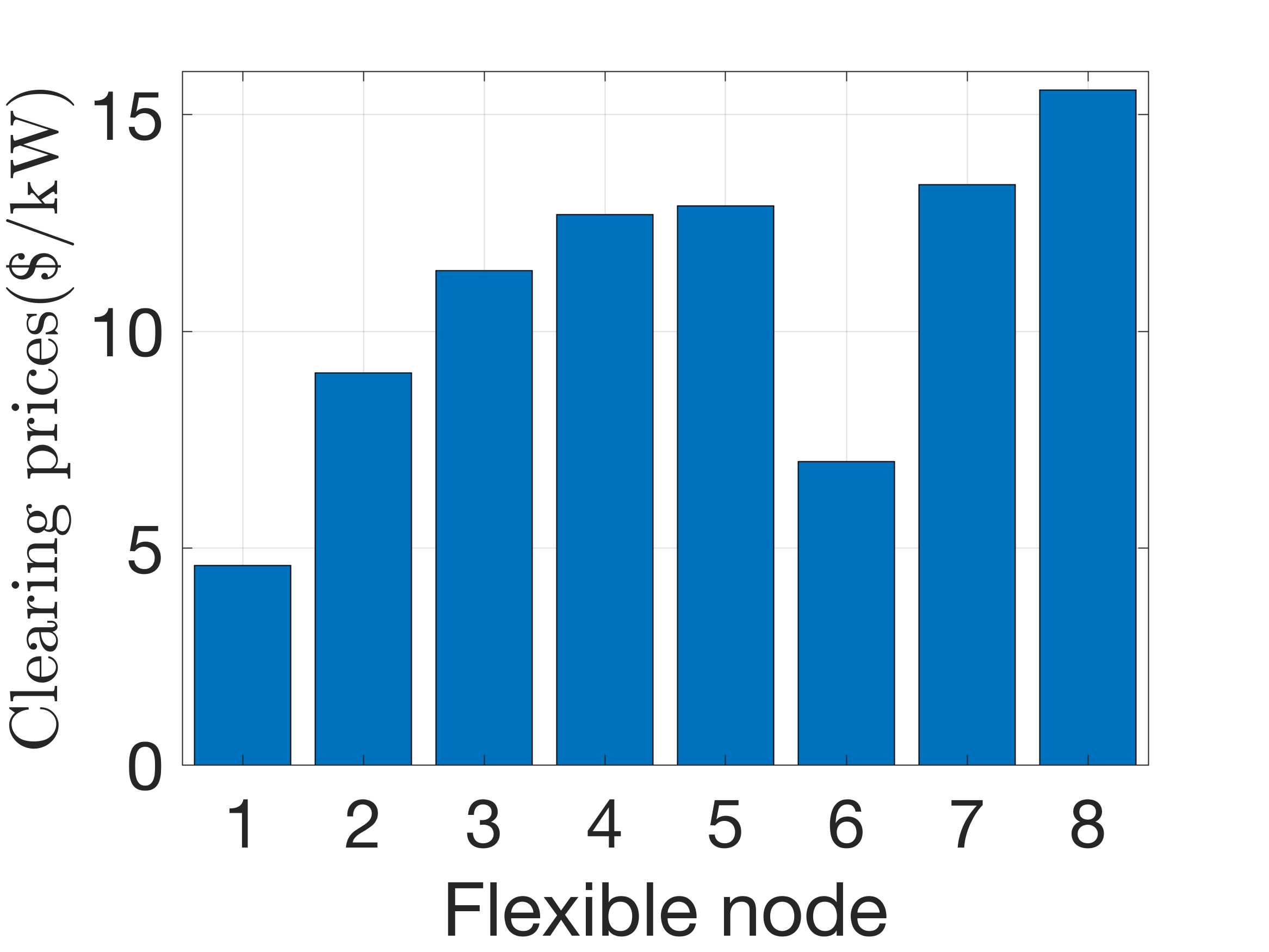}\label{fig:clearing_uncer}}
\caption{(a) Allocation of Aggregator flexibility after market clearing under varying Aggregator prices and demand uncertainty, (b) Clearing price based on cheapest allocation at each node under uncertainty in demand.}
\end{figure}

In the previous sections, the simulations assumed a known background demand at each node.
However, in forecasting the background demand at each node, we need to account for prediction errors to ensure a grid-aware allocation of flexibility. 
In this section, we consider a robust formulation by accounting for the uncertainty at each flexible node as shown in Table~\ref{table_unc}. Based on this demand uncertainty, we solve (ROB-P1) and (ROB-P2) to obtain the grid-aware Aggregator capacity allocation, which is shown in Fig.~\ref{fig:Agg_alloted_uncer}. From Fig.~\ref{fig:Agg_alloted_uncer}, the total allocation is reduced by $0.44$ MW in this case (compared to the deterministic case) as (ROB-P1) and (ROB-P2) need to reserve certain amount of network capacity to account for the uncertainty in background demand. Based on this allocation, we again determine the clearing price based on the least price Aggregator at each flexible node (which is again higher than the price in the deterministic case) and this clearing price for the robust case is depicted in Fig.~\ref{fig:clearing_uncer}. Since robustness raises the cleared price but reduces flexible capacity, there is a trade-off when calculating DNO revenue. In this case the DNO revenue is calculated as $\$39,540$, which is 4.8\% less than the deterministic case. \textcolor{black}{Even though the simulation results in this paper utilized a modified IEEE-37 node system, the CIA formulation is convex and has been shown to scale to larger networks~\cite{nazir2020gridaware}.}

\begin{table}[!t]
\caption{Uncertainty in background demand at flexible nodes (MW)}
\label{table_unc}
\centering
\begin{tabular}{rllllllll}
\hline
{Node} & {1} & {2} & {3} & {4} & {5} & {6} & {7} & {8} \\
\hline
Uncertainty & 0.2 & 0.1 & 0.02 & 0.1 & 0.2 & 0.03 & 0.02 & 0.03 \\
\hline
\end{tabular}
\end{table}

\section{Conclusions and future work}
 Leveraging recent work on convex inner approximation, we have now formulated a novel grid-aware allocation market-clearing mechanism which is executed in two steps. Step~1 considers a feasibility check to determine whether the available Aggregator flexibility can be safely hosted on the network. In Step~2, flexibility is prioritized based on Aggregator's bids and AC network physics and constraints to provide grid-aware allocation across nodes in the network. 

Future work will extend the work on uncertainty in background (net) demand, especially at non-flexible nodes, to consider DNO margins and Aggregator flexibility reserves. Recent work by~\cite{Dongchan} offers an interesting and systematic approach to develop a robust convex inner approximation that can be employed in (P1) and (P2). Finally, the authors are pursuing a hybrid solution to manage uncertainty that considers that the Utility can share real-time grid measurements with the Aggregators, which can then adapt its cleared grid-aware allocations dynamically. This approach raises the prospect that Aggregators can then offer the DNO flexibility reserves to incentivize improved forecasts of background demand.   



%



 \bibliographystyle{IEEEtran}
 \small\bibliography{sample.bib}

\end{document}